\documentclass[letterpaper, 10 pt, conference]{ieeeconf}  

\IEEEoverridecommandlockouts                              
\usepackage{graphicx}

 \usepackage[pdftex]{color}
\usepackage{amsmath}
\usepackage{amsfonts} 
\usepackage{algorithmic}
\usepackage{array}
\usepackage{soul}
\newtheorem{teor}{Theorem}[section]
\newtheorem{corr}{Corollary}[section]
\newtheorem{propo}{Proposition}[section]
\newtheorem{lemm}{Lemma}[section]

\newcommand{\tp}{^{\top}}

\newcommand{\beq}{\begin{equation}}
\newcommand{\eeq}{\end{equation}}
\newcommand{\bea}{\begin{eqnarray}}
\newcommand{\eea}{\end{eqnarray}}

\newcommand{\bsea}{\begin{subeqnarray}}
\newcommand{\esea}{\end{subeqnarray}}
\newcommand{\nn}{\nonumber}

\def\bmat{\left[ \begin{array}}
\def\emat{\end{array} \right]}
\DeclareMathOperator{\tr}{tr} 
\DeclareMathOperator{\diag}{diag} 
\definecolor{Royalblue}{cmyk}{1,0.30,0.2,0.2}

\newcommand{\alg}[1]{\begin{align}#1\end{align}}

\begin{document}

\title{\Huge{Factor analysis with finite data }}

\author{Valentina Ciccone, Augusto Ferrante, and Mattia Zorzi\thanks{V. Ciccone, A. Ferrante, and M. Zorzi are with   the Department of Information
Engineering, University of Padova, Padova, Italy; e-mail: {\tt\small valentina.ciccone@dei.unipd.it} (V. Ciccone); {\tt\small augusto@dei.unipd.it} (A. Ferrante); {\tt\small zorzimat@dei.unipd.it} (M. Zorzi).}}


\maketitle

\thispagestyle{empty}
\pagestyle{empty}


\begin{abstract}

Factor analysis aims to describe high dimensional random vectors by means of a small number of unknown common factors. In mathematical terms, it is required to decompose the covariance matrix $\Sigma$ of the random vector  as the sum of a diagonal matrix $D$ | accounting for the idiosyncratic noise in the data | and a low rank matrix $R$ | accounting for the variance of the common factors | in such a way that the rank of $R$ is as small as possible so that the number of common factors is minimal.

In practice, however, the matrix $\Sigma$ is unknown and must be replaced by its estimate, i.e. the sample covariance, which comes from a finite amount of data. This paper provides a strategy to account for the uncertainty in the estimation of $\Sigma$ in the factor analysis problem. 
\end{abstract}

\section{Introduction}

Factor models are used to summarize high-dimensional data vectors with a small number of 
unknown and non-observed common factors. They boast a long tradition in different disciplines such as psychometrics, econometrics, systems identification and control engineering. 
The history of these models can be tracked back to the beginning of the last century in the framework of  psychological tests \cite{spearman_1904,BURT_1909,THURSTONE_1934} and, since then, their importance has spread in virtually all disciplines of sciences \cite{KALMAN_SELECTION_ECONOMETRICS_1983,SCHUPPEN_1986,Bekker-deLeeuw_1987,PICCI_1987,GEWEKE_DYNAMIC,PICCI_PINZONI_1986,Pena_BOX__1987,ON_THE_PENA_BOX_MODEL_2001,deistler2007,DEISTLER_1997,ANDERSON_DEISTLER_1993,SARGENT_SIM_1977,forni_lippi_2001,ENGLE_ONE_FACTOR_1981,WATSON_ALGORITHMS_1983,MFA}; see also the more recent papers \cite{LATENTG,BSL,Bottegal-Picci,DEISTLER_2015} where a larger number of other references are listed. 
Furthermore, the mathematical analysis of these models has been carried out by several different perspectives: for example, a detailed geometric description of such models can be found in \cite{scherrer1998structure}, while a maximum likelihood approach in a statistical testing framework has been proposed for example in the seminal work \cite{anderson1956statistical}.

\begin{figure}[htb]
	\centering
	\includegraphics[width=8cm]{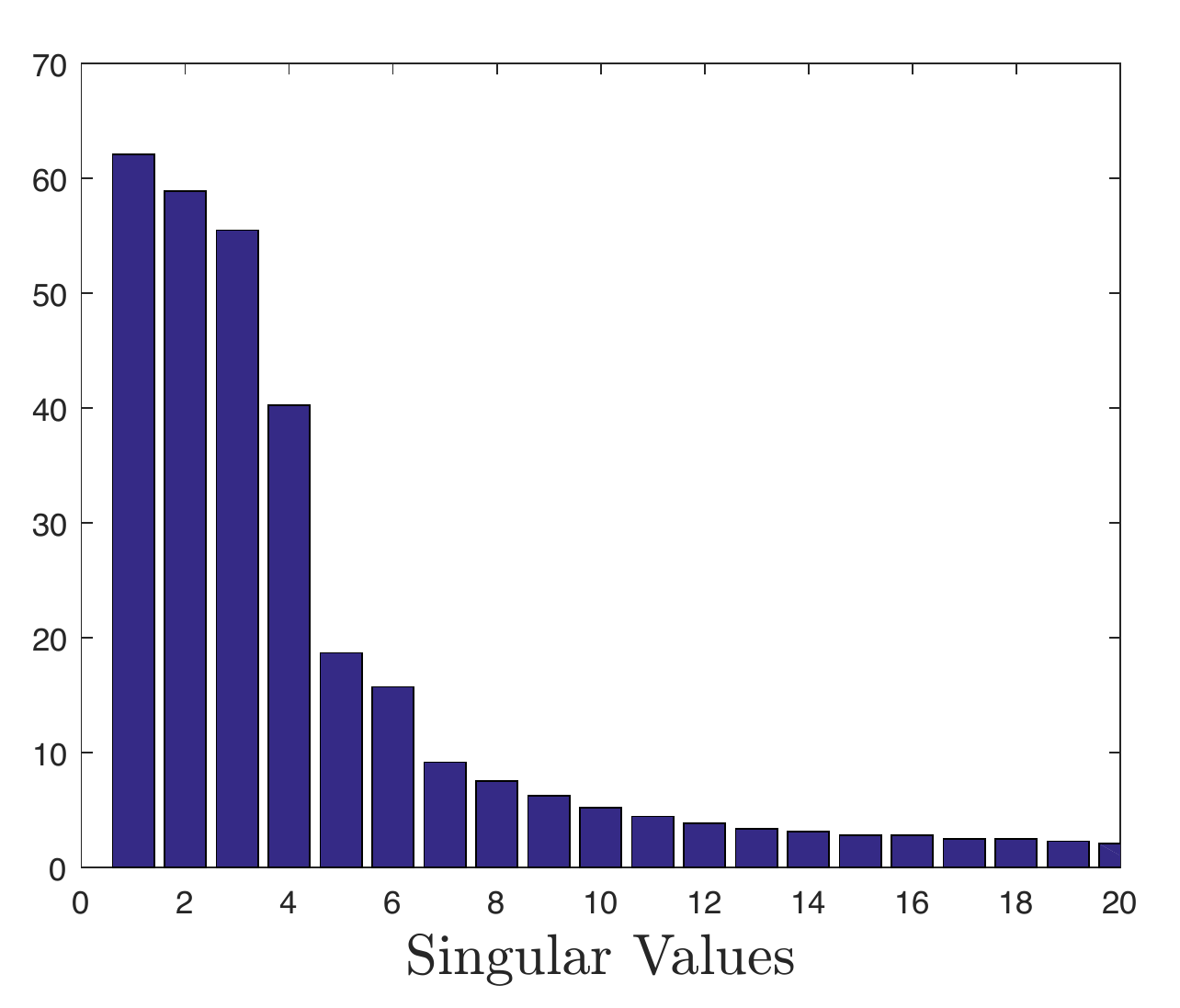}
	\caption{First $20$ singular values of the low-rank matrix $R$ obtained by applying the standard factor analysis decomposition algorithm to the estimated covariance matrix $\hat{\Sigma}$.}
\label{fig:Figone}
\end{figure}

In its classical and most simple static version the construction of a factor model can be reduced to a particular case of high-dimensional matrix additive decomposition problem. This kind of problems arise naturally in numerous frameworks and have therefore received a great deal of attention,
see \cite{Chandrasekaran-Sanghavi-Parrilo-Willsky,agarwal2012noisy,Zorzi-Ferrante-Automatica-12,Ferrante-Pav-Zorzi-12} and references therein.
More precisely, for the identification of a factor model we assume that the covariance matrix $\Sigma$ of a high-dimensional vector of data is assigned and must be decomposed as $\Sigma=R+D$, that is the sum of a non-negative diagonal matrix $D$ modelling the covariance of the idiosyncratic noise and a positive semidefinite matrix $R$, having rank as small as possible, modelling the covariance of the latent variable which turns out to be a combination of the common factors.
For this mathematical problem to be meaningful with respect to the original factor model, we have to assume that either $\Sigma$ is known with good degree of precision or that small perturbations of $\Sigma$ have no or little effect on the rank of the corresponding matrix $R$ in the additive decomposition of $\Sigma$.
It seems fair to say that both of these conditions are usually not satisfied:
$\Sigma$ must be estimated from the available data that are certainly finite and possibly limited so that we can expect an estimation error whose covariance can be usually estimated with reasonable precision.
On the other hand, simulations show that the rank of $R$ is rather sensitive to pertubations of $\Sigma$. For example we have considered a model with 1000 samples of a 50-dimensional vector of data generated by $4$ non-observed common factors.
By applying the standard factor analysis decomposition algorithm (based on the minimization of the nuclear norm) to the estimated covariance matrix $\hat{\Sigma}$, we obtained a matrix $R$ whose first $20$ singular values are depicted in 
Figure \ref{fig:Figone}.
It is clear that this matrix is far from having rank $4$ as it would be if the procedure returned the correct model.
For a sanity check, we have also implemented the standard factor analysis decomposition algorithm to the true covariance matrix $\Sigma$ and obtained a matrix $R$ in which the fifth singular value is 
$10^6$ times smaller than the fourth.

Motivated by this issue, we have considered the problem of taking the uncertainty in the estimation of $\Sigma$ into account.
This leads to a much more complex problem where we are required to {\em compute} the matrix $\Sigma$
in such a way that the rank of $R$ in the additive decomposition $\Sigma=R+D$ is minimized under a constraint limiting the Kullback-Leibler divergence between $\Sigma$ and the estimated covariance $\hat{\Sigma}$ to a prescribed tolerance that depends on the precision of our estimate  $\hat{\Sigma}$.
The problem in this formulation appears to be quite hard to solve as the number of variables is large the constraints are difficult to impose as the solution will always lie on the boundary.
Our strategy is to resort to the dual analysis which is rather delicate to carry over but 
yields a problem that appears much easier to tackle. Moreover it provides a necessary and sufficient condition for the uniqueness of the solution of the original problem.

The paper is organized as follows. In the Section \ref{sec_pb_form} we recall the classical factor analysis problem and, from it, we derive a mathematical formulation of our generalized factor analysis  problem. In Section \ref{sec_dual} we derive a dual formulation of our problem. In Section \ref{sec_ex} we prove existence and uniqueness of the solution for the dual problem. Then, in Section \ref{sec_recovery} we recover the solution of the primal problem. Finally, some conclusions are provided.
We warn the reader that the present paper only reports some preliminary result. In particular, all the proofs are omitted and will appear in a forthcoming and more complete publication.

{\em Notation:} Given a vector space $\cal V$ and a subspace ${\cal W}\subset {\cal V}$, we denote by
${\cal W}^\bot$ the orthogonal complement of ${\cal W}$ in ${\cal V}$.
Given a matrix $M$, we denote its transpose by $M\tp$; if $M$ is a square matrix 
$\tr(M)$ denotes its trace, i.e. the sum of the elements in the main diagonal of $M$;
moreover, $|M|$ denotes the determinant of $M$. {  We denote the spectral norm of $M$ as $\Vert M \Vert_2$.}
We endow the space of square real matrices with the following inner product: for $A,B\in{\mathbb R}^{n\times n}$, $\langle A,B\rangle:=\tr(AB)$. 
The kernel of a matrix (or of a linear operator) is denoted by $\ker(\cdot)$.
The symbol $\mathbf{Q}_n$ denotes the vector space of real symmetric matrices of size $n$. If $X\in\mathbf{Q}_n$ is positive definite or positive semi-definite we write $X\succ 0$ or $X\succeq 0$, respectively.
Moreover, we denote by $\mathbf{D}_n$ the vector space of  diagonal matrices of size $n$; $\mathbf{D}_n$ is clearly a subspace of $\mathbf{Q}_n$ and we denote by
$\mathbf{M}_n:=\mathbf{D}_n^{\bot}$ the orthogonal complement of $\mathbf{D}_n$ in $\mathbf{Q}_n$ (with respect to the inner product just defined).
It is easy to see that 
$\mathbf{M}_n$ is the vector space of symmetric matrices of size $n$ with only zero elements on the main diagonal.\\
Moreover, we denote by $\diag(\cdot)$ the operator mapping $n$ real elements $d_i, i=1,...,n$ into the diagonal matrix having the $d_i$'s as elements in its main diagonal and, given a linear operator $\chi(\cdot)$, we denote by $\chi^*(\cdot)$ the corresponding adjoint operator. 

\section{Problem Formulation}\label{sec_pb_form}

We consider a standard,  static linear  factor model that can be represented as follows:
 \alg{x&=Aw_y+Bw_z\nn\\
 y&=A w_y\nn\\
 z&=Bw_z\nn}
 where $A\in \mathbb{R}^{n\times r}$ with $r<<n, B \in \mathbb{R}^{n\times n}$ diagonal. $A$ is the factor loading matrix, $y$ is the latent variable, and $Bw_z$ is the idiosyncratic noise component. $w_y$ and $w_z$  are independent Gaussian random vectors with zero mean and covariance matrix equal to the identity matrix of dimension $r$ and $n$, respectively. Note that, $w_y$ represents the (independent) latent factors. 
Consequently, $x$ is a Gaussian random vector with zero mean; we denote by $\Sigma$ its covariance matrix. Since $y$ and $z$ are independent we get that
\begin{equation}
\label{decomposition}
\Sigma=R+D
\end{equation}
where $R$ and $D$ are the covariance matrices of $y$ and $z$, respectively. Thus, $R:=AA\tp $  has rank equal to $r$, and $D=BB\tp $ is diagonal.\\
The objective of factor analysis consists in finding a decomposition ``low-rank plus diagonal'' \eqref{decomposition} of $\Sigma$. This amounts to solve the minimum rank problem
\begin{equation}
\label{rank_min}
\begin{aligned}
\min_{R,D\in \mathbf{Q}_n} \quad & \text{ rank}(R)\\
\text{subject to }\quad & R,D\succeq 0\\
					& D \in \mathbf{D}_n\\
					& {\Sigma}=R+D
\end{aligned}
\end{equation}
which is, however, an hard problem. 
A well-known convex relaxation of \eqref{rank_min} is the trace minimization problem
 \begin{equation}
\label{original_problem}
\begin{aligned}
\min_{R,D\in \mathbf{Q}_n} \quad & \text{ tr}(R)\\
\text{subject to }\quad & R,D\succeq 0\\
					    & D\in \mathbf{D}_n\\
					    & {\Sigma}=R+D.
\end{aligned}
\end{equation}
Problem (\ref{original_problem}) provides a solution which is a good approximation of the one of Problem (\ref{rank_min}), \cite{MIN_RANK_SHAPIRO_1982}. This evidence is justified by the fact that $\tr(R)$, i.e. the nuclear norm of $R$, is the convex hull of $\text{ rank}(R)$ over the set ${\cal S}:=\{ R\in \mathbf{Q}_n \; \mathrm{s.t.} \; \| R\|_2\leq 1\}$, \cite{FAZEL_MIN_RANK_APPLICATIONS_2002}.

In practice, however, the matrix $\Sigma$ is not known and needs to be estimated from a $N$-length realization (i.e. a data record) $\mathrm{x}_1 \ldots \mathrm{x}_N$ of $x$. The typical choice is to take the sample covariance estimator \alg{\hat \Sigma:=\frac{1}{N}\sum_{k=1}^N  \mathrm{x}_k  \mathrm{x}_k\tp }
which is statistically consistent, i.e. $\hat \Sigma$ almost surely converges to $\Sigma$ as $N$ tends to infinity. If we replace $\Sigma$ with $\hat \Sigma$ in (\ref{original_problem}) then, 
as discussed in the Introduction the corresponding solution will rapidly degrade unless we provide an appropriate model accounting for the error in the estimation of $\Sigma$.
Let $\hat x$ be a Gaussian random vector with zero mean and covariance matrix $\hat \Sigma$. Note that there exists a one to one correspondence between $x$ and $\Sigma$, and between $\hat x$ and $\hat \Sigma$; $\hat x$ is a crude ``model approximation'' for $x$. Thus, to account for this uncertainty, we assume that  $\Sigma$ belongs to a ``ball'' of radius $\delta /2$ centered in $\hat \Sigma$. Such a ball is formed by placing a bound on the the {\em Kullback-Leibler} divergence between $x$ and $\hat x$:
\alg{{\cal B}:=\{ \Sigma\in \mathbf{Q}_n \; \mathrm{s.t.}\; \Sigma\succ 0,\;{\cal D}_{KL}(\Sigma\|\hat \Sigma)\leq \delta /2\}.} 
Here ${\cal D}_{KL}$ is the {\em Kullback-Leibler} divergence defined by:
\alg{ {\cal D}_{KL}(\Sigma\|\hat \Sigma):= \frac{1}{2}\left(-\log |\Sigma|+\log |\hat \Sigma|+\tr(\Sigma \hat \Sigma^{-1})-n\right).\nn} This way to deal with model uncertainty has been successfully in econometrics for model mispecification \cite{ROBUSTNESS_HANSENSARGENT_2008} and in robust filtering \cite{ROBUST_FILTERING,ROBUST_TAU,ROBUST_KALMAN_2017,ROBUST_CONV_2015,convtau}. Accordingly, the trace minimization problem can be reformulated as follows 
\begin{equation}
\label{primal_before_changing_variables}
\begin{aligned}
\min_{\Sigma,R,D\in  \mathbf{Q}_n} & \text{ tr}(R)\\
\text{subject to  } &  R,D\succeq 0\\
					& D\in \mathbf{D}_n\\
					& \Sigma=R+D\\
					& \Sigma\in \cal B.
\end{aligned}
\end{equation}

Note that, in (\ref{primal_before_changing_variables}) we can eliminate variable $D$, obtaining the equivalent problem 
\begin{equation}
\label{primal}
\begin{aligned}
\min_{R,\Sigma\in  \mathbf{Q}_n} \quad & \text{ tr}(R)\\
\text{subject to  } \quad & R,\Sigma-R\succeq 0 \\
					& \chi(\Sigma-R)=0 \\
					& \Sigma\succ 0\\
					&2 \mathcal{D}_{KL}(\Sigma||\hat{\Sigma})\leq \delta  
\end{aligned}
\end{equation}
where $\chi(\cdot)$ is the self-adjoint operator orthogonally projecting onto $\mathbf{M}_n$, i.e. if $M\in  \mathbf{Q}_n$, $\chi(M)$ is the matrix of $\mathbf{M}_n$ in which each off-diagonal element is equal to the corresponding element of $M$ (each diagonal element of $\chi(M)$ is clearly zero).

\subsection{The choice of $\delta$}
The tolerance $\delta$ may be chosen by taking into account the accuracy of the estimate of $\Sigma$ which, in turn, depends on the numerosity of the available data.  
Notice, however, that if $\delta$ in \eqref{primal} is sufficiently large, we may obtain an optimal solution such that  R is identically equal to the null matrix and $\Sigma\in \mathbf{D}_n$.
In order to avoid this trivial situation we need to require that the maximum tolerable  Kullback-Leibler divergence $\delta$ in \eqref{primal} is strictly less than a certain $\delta_{max}$ that can be determined as follows: since the trivial solution $R= 0$ would imply a diagonal $\Sigma$, that is $\Sigma=\Sigma_D:=\text{diag}(d_1,...,d_n)>0$, $\delta_{max}$ can be determined by solving the following minimization problem
\begin{align}
\label{min_diagonal}
\delta_{max}:=\min_{\Sigma_D\in\mathbf{D}_n  } \mathcal{D}_{KL}(\Sigma_D \Vert \hat{\Sigma}).\end{align}
\begin{propo}\label{prop_deltamax} Let $\gamma_i$ denote the element $i$-th element in the main diagonal of the inverse of the sample covariance $\hat{\Sigma}^{-1}$. Then, the optimal $\Sigma_D$ which solves the minimization problem in \eqref{min_diagonal} is given by
\[\Sigma_D= \text{diag}(\gamma_1^{-1}, ...,\gamma_n^{-1}).\]
Moreover, $\delta_{max}$ can be determined as 
\begin{equation}
\label{delta_MAX}
\delta_{max}=\mathcal{D}_{KL}(\Sigma_D \Vert \hat{\Sigma})= \log|[\hat{\Sigma}^{-1}-\chi(\hat{\Sigma}^{-1})]\hat{\Sigma}|.
\end{equation} \end{propo} In what follows, we always assume that 
$\delta$ in \eqref{primal} strictly less than $\delta_{max}$, so that the trivial solution $R\equiv 0$ is ruled out.
\section{Dual Problem}\label{sec_dual}
We start by formulating the constrained optimization problem in  \eqref{primal} as an unconstrained minimization problem. 
The Lagrangian associated to \eqref{primal} is
\begin{equation}
\label{lagrangian}
\begin{aligned}
\mathcal{L} (R,& \Sigma, \lambda, \Lambda, \Gamma, \Theta ) \\
= & \text{tr}(R) +\lambda(-\log|\Sigma|+\log|\hat{\Sigma}|-n +\text{tr}(\hat{\Sigma}^{-1}\Sigma)-\delta)\\
& -\text{tr}(\Lambda R)-\text{tr}(\Gamma(\Sigma-R))+\text{tr}{(\Theta \chi(\Sigma-R)})\\
= & \text{tr}(R) +\lambda(-\log|\Sigma|+\log|\hat{\Sigma}|-n +\text{tr}(\hat{\Sigma}^{-1}\Sigma)-\delta)\\
& -\text{tr}(\Lambda R)-\text{tr}(\Gamma(\Sigma-R))+\text{tr}{(\chi^{*}(\Theta)(\Sigma-R)})\\
= & \text{tr}(R) +\lambda(-\log|\Sigma|+\log|\hat{\Sigma}|-n +\text{tr}(\hat{\Sigma}^{-1}\Sigma)-\delta)\\
& -\text{tr}(\Lambda R)-\text{tr}(\Gamma(\Sigma-R))+\text{tr}{(\chi(\Theta)(\Sigma-R)})\\
\end{aligned}
\end{equation}
with $\lambda\in\mathbb{R}, \lambda\geq 0$, and $\Lambda,\Gamma, \Theta \in \mathbf{Q}_n$ with $\Lambda,\Gamma\succeq 0$. 
\\ 
The first and the second equality are due to the fact that the operator $\chi(\cdot)$ is self-adjoint.\\
Notice that in \eqref{lagrangian} we have not included the constraint $\Sigma\succ 0$. This is due to the fact that, as we will see later on, such condition is automatically fulfilled by the solution of the dual problem.\\
The dual function is now the infimum of $\mathcal{L}(R,\Sigma, \lambda, \Lambda, \Gamma, \Theta )$ over $R$ and $\Sigma$.

Since the Langrangrian is convex in order to find the minimum we use standard variational methods.\\
The first variation of the Lagrangian \eqref{lagrangian} at $\Sigma$ in direction $\delta\Sigma\in\mathbf{Q}_n$ is 
\begin{equation*}
\delta\mathcal{L}(\Sigma;\delta\Sigma) = \text{tr}(-\lambda\Sigma^{-1}\delta\Sigma+\lambda\hat{\Sigma}^{-1}\delta\Sigma-\Gamma\delta\Sigma+\chi(\Theta)\delta\Sigma).
\end{equation*}
By imposing the optimality condition 
\[\delta\mathcal{L}(\Sigma;\delta\Sigma)=0, \qquad\forall\delta\Sigma\in\mathbf{Q}_n,\]
 which is equivalent to impose $\text{tr}(-\lambda\Sigma^{-1}\delta\Sigma+\lambda\hat{\Sigma}^{-1}\delta\Sigma-\Gamma\delta\Sigma+\chi(\Theta)\delta\Sigma) =0$ for all $\delta\Sigma\in\mathbf{Q}_n$, we obtain 
\begin{equation}
\label{der_var_wrt_sigma_opt}
\begin{aligned} 
\Sigma=\lambda(\lambda\hat{\Sigma}^{-1}-\Gamma+\chi(\Theta))&^{-1}
\end{aligned}
\end{equation}
provided that $\lambda\hat{\Sigma}^{-1}-\Gamma+\chi(\Theta)\succ 0$ and $\lambda>0$. Note that these conditions are equivalent to impose that the optimal $\Sigma$ that minimizes the Lagrangian satisfies the constraint $\Sigma\succ 0$.\\

The first variation of the Lagrangian \eqref{lagrangian} at $R$ in direction $\delta R\in\mathbf{Q}_n$ is
\begin{equation*}
\delta\mathcal{L}(R;\delta R)= \text{tr}(\delta R - \Lambda\delta R+\Gamma\delta R -\chi(\Theta)\delta R).
\end{equation*}
Again, by imposing the optimality condition
\[\delta\mathcal{L}(R;\delta R)=0, \qquad \forall\delta R\in\mathbf{Q}_n,\]
which is equivalent to $\text{tr}(\delta R - \Lambda\delta R+\Gamma\delta R -\chi(\Theta)\delta R)= 0$ for all $\delta R\in\mathbf{Q}_n$, we get that 
 \begin{equation}
\label{der_var_wrt_r}
\begin{aligned}
I -\Lambda +\Gamma -\chi(\Theta)&=0 .
\end{aligned}
\end{equation}

\begin{propo} \label{prop_dual}The dual problem of (\ref{primal}) is 
\begin{equation}
\label{max_dual}
\max_{(\lambda,\Gamma, \Theta) \in \mathcal{C}_0} J(\lambda,\Gamma,\Theta) 
\end{equation}
where 
\begin{equation*}
\label{dual}
\begin{aligned}
J(\lambda,& \Gamma,\Theta):= \lambda(\log|(\hat{\Sigma}^{-1}+\lambda^{-1}(\chi(\Theta)-\Gamma))|+\log|\hat{\Sigma}|-\delta)
\end{aligned}
\end{equation*}
and $\mathcal{C}_0$ is defined as 
\alg{
\label{C0}
\mathcal{C}_0:=\lbrace (\lambda,\Gamma,\Theta):\ &\lambda>0, \ I+\Gamma-\chi(\Theta)\succeq0, \   
\Gamma\succeq 0,\nn\\ &\hat{\Sigma}^{-1}+\lambda^{-1}(\chi(\Theta)-\Gamma)\succ 0\rbrace.
}
\end{propo}
Note that, the conditions $\lambda>0$ and $\hat{\Sigma}^{-1}+\lambda^{-1}(\chi(\Theta)-\Gamma)\succ 0$ arise from \eqref{der_var_wrt_sigma_opt}.

\section{Existence and uniqueness of the solution for the dual problem} \label{sec_ex}

We reformulate the maximization problem in \eqref{max_dual} as a minimization problem.
\begin{equation}
\label{Dual}
\begin{aligned}
\min_{(\lambda,\Gamma, \Theta) \in \mathcal{C}_0} \tilde{J}(\lambda,\Gamma, \Theta )
\end{aligned}
\end{equation}
where 
\alg{\tilde J(\lambda,\Gamma, \Theta)& = \lambda(-\log|\hat{\Sigma}^{-1}+\lambda^{-1}(\chi(\Theta)-\Gamma)|-\log|\hat{\Sigma}|+\delta). \nn}

\subsection{Existence}
As it is often the case, existence of the optimal solution is a very delicate issue.
Our strategy in order to deal with this issue and prove that the dual problem \eqref{Dual} admits a solution consists in showing that we can restrict our set $\mathcal{C}_0$ to a smaller compact set $\mathcal{C}$ over which the minimization problem is equivalent to the one in \eqref{Dual}. Since the objective function is continuous over $\mathcal{C}_0$, and hence over $\mathcal{C}$, by Weirstrass's theorem $\tilde{J}$ admits a minimum.

First, we recall that the operator $\chi(\cdot)$ is self-adjoint. Moreover, we notice that $\chi(\cdot)$ is not injective on $\Theta$, thus we want to restrict the domain of $\chi(\cdot)$ to those $\Theta$ such that $\chi(\cdot)$ is injective. Since $\chi$ is self-adjoint we have that:
\[\ker(\chi)=[\text{range }\chi]^{\perp}.\]
Thus, by restricting $\Theta$ to range($\chi$)$=[\ker(\chi)]^{\perp}=\mathbf{M}_n$, the map becomes injective. Therefore, without loss of generality, from now on we can safely assume that $\Theta\in\mathbf{M}_n$ so that
$\chi(\Theta)=\Theta$ and we 
restrict our set $\mathcal{C}_0$ to $\mathcal{C}_{1}$:
\begin{align*}
\mathcal{C}_{1}: =& \lbrace (\lambda,\Gamma,\Theta)\in \mathcal{C}_0: \Theta \in \mathbf{M}_n \rbrace \\
= & \lbrace (\lambda,\Gamma,\Theta):\lambda>0, I+\Gamma-\Theta\succeq 0, \Gamma\succeq 0, \\
& \Theta \in \mathbf{M}_n,\,  \hat{\Sigma}^{-1}+\lambda^{-1}(\Theta-\Gamma)\succ 0 \rbrace.
\end{align*}
Moreover, since $\Theta$ and $\Gamma$ enter into the problem always through their difference they cannot be univocally determined individually. However, their difference does. This allows us to restrict $\Gamma$ to the space of the diagonal positive semi-definite matrices.
For this reason, we can further restrict our set $\mathcal{C}_1$ to $\mathcal{C}_{2}$:
\begin{align*}
\mathcal{C}_{2}:= & \lbrace (\lambda, \Gamma, \Theta): \lambda>0, I+\Gamma-\Theta\succeq 0, \Gamma\succeq 0, \Gamma\in\mathbf{D}_n,\\
& \Theta \in \mathbf{M}_n, \hat{\Sigma}^{-1}+\lambda^{-1}(\Theta-\Gamma)\succ 0 \rbrace .
\end{align*}

To further restrict this set, we need to find a lower bound on $\lambda$ which has an infimum but not a minimum on $\mathcal{C}_{2}$.
The following result provides such a bound.
\begin{lemm}\label{lem41}
Let $(\lambda_k, \Gamma_k, \Theta_k )_{k\in \mathbb{N}}$ be a sequence of elements in $\mathcal{C}_2$ such that
$$
\lim_{k\rightarrow \infty}
\lambda_k = 0.$$
Then $(\lambda_k, \Gamma_k, \Theta_k )_{k\in \mathbb{N}}$ is not an infimizing sequence for $\tilde{J}$ .
\end{lemm}

As a consequence of the previous result, we have that minimization of the dual functional over the set $\mathcal{C}_2$ is equivalent to minimization over the  set:
\begin{align*}
\mathcal{C}_{3}:= &\lbrace (\lambda,\Gamma,\Theta):\lambda\geq \varepsilon , I+\Gamma-\Theta\succeq 0, \Gamma\succeq 0, \Gamma\in\mathbf{D}_n,\\
& \Theta \in \mathbf{M}_n , \hat{\Sigma}^{-1}+\lambda^{-1}(\Theta-\Gamma)\succ 0 \rbrace &
\end{align*}
for a certain $\varepsilon>0$.\\

\noindent The next result is a counterpart of the previous one as it deals with the fact that, so far, $\lambda$ is still unbounded and thus there could, in principle, exist an infimizing sequence for which the corresponding  $\lambda$ diverges. This is not the case in view of the following lemma.


\begin{lemm}
Let $(\lambda_k, \Gamma_k, \Theta_k )_{k\in \mathbb{N}}$ be a sequence of elements in $\mathcal{C}_3$ such that
\beq
\label{ltinf}
\lim_{k\rightarrow \infty}
\lambda_k= \infty.
\eeq
Then $(\lambda_k, \Gamma_k, \Theta_k )_{k\in \mathbb{N}}$ is not an infimizing sequence for $\tilde{J}$ . 
\end{lemm}

\bigskip

As a consequence, the feasible set $\mathcal{C}_3$ can be further restricted to the set:
\begin{align*}
\mathcal{C}_{4}:=&  \lbrace (\lambda, \Gamma, \Theta):\varepsilon \leq \lambda\leq M , I+\Gamma-\Theta\succeq 0, \Gamma\succeq 0,\\ 
& \Gamma\in\mathbf{D}_n, \Theta \in \mathbf{M}_n , \hat{\Sigma}^{-1}+\lambda^{-1}(\Theta-\Gamma)\succ 0 \rbrace 
\end{align*}
for a certain $M<\infty$.\\

\noindent The next result provides an upper bound for $\Theta-\Gamma$.


\begin{lemm}
Let $(\lambda_k, \Gamma_k, \Theta_k )_{k\in \mathbb{N}}$ be a sequence of elements in $\mathcal{C}_4$ such that
\beq
\label{tmgtinf}
\lim_{k\rightarrow \infty}
\Vert \Theta_k-\Gamma_k\Vert = +\infty.
\eeq
Then $(\lambda_k, \Gamma_k, \Theta_k )_{k\in \mathbb{N}}$ is not an infimizing sequence for $\tilde{J}$. 
\end{lemm}

It follows that  there exists $\rho$ s.t. $|\rho|<\infty$ and 
\[ \Theta -\Gamma\succeq \rho I.\]
Therefore, the feasible set $\mathcal{C}_4$ can be further restricted to become:\\
\begin{align*}
\mathcal{C}_{5}:= & \lbrace (\lambda,\Gamma,\Theta): \varepsilon \leq \lambda\leq M , \rho I\preceq \Theta -\Gamma\preceq I, \Gamma\succeq 0, \\
& \Gamma\in\mathbf{D}_n, \Theta \in \mathbf{M}_n , \,\hat{\Sigma}^{-1}+\lambda^{-1}(\Theta-\Gamma)\succ 0 \rbrace. 
\end{align*}

Now observe that in $\mathcal{C}_{5}$ $\Theta$ and $\Gamma$ are orthogonal so that 
if $(\lambda_k, \Gamma_k, \Theta_k )_{k\in \mathbb{N}}$ is a sequence of elements in $\mathcal{C}_5$ such that
\beq
\label{gtinf}
\lim_{k\rightarrow \infty}
\Vert \Gamma_k \Vert = +\infty
\eeq
or
\beq
\label{ttinf}
\lim_{k\rightarrow \infty}
\Vert \Theta_k \Vert = +\infty
\eeq
then (\ref{tmgtinf}) holds.
Then we have the following

\begin{corr}
Let $(\lambda_k, \Gamma_k, \Theta_k )_{k\in \mathbb{N}}$ be a sequence of elements in $\mathcal{C}_5$ such that (\ref{gtinf}) or (\ref{ttinf}) holds. Then $(\lambda_k, \Gamma_k, \Theta_k )_{k\in \mathbb{N}}$ is not an infimizing sequence for $\tilde{J}$.
\end{corr}

 
Thus minimizing over the set $\mathcal{C}_{5}$ is equivalent to minimize over:
\begin{align*}
\mathcal{C}_{6}:= & \lbrace (\lambda, \Gamma, \Theta): \varepsilon \leq\lambda\leq M , \rho I\preceq \Theta -\Gamma\preceq I, 0\preceq\Gamma\preceq \alpha I, \\
& \Gamma\in\mathbf{D}_n,\Theta \in \mathbf{M}_n , \, \hat{\Sigma}^{-1}+\lambda^{-1}(\Theta-\Gamma)\succ 0 \rbrace 
\end{align*}
for a certain  $\alpha$  such that $0<\alpha<+\infty$.\\

Finally, let us consider a sequence $(\lambda_k, \Gamma_k, \Theta_k )_{k\in \mathbb{N}}\in\mathcal{C}_6$ such that, as $k\rightarrow\infty$, $\big{|}\hat{\Sigma}^{-1}+\lambda_k^{-1}(\Theta_k-\Gamma_k)\big{|}\rightarrow 0$. This implies that $\tilde{J}\rightarrow +\infty$. Thus, such sequence does not infimize the dual functional.
Thus, the final feasible set $\mathcal{C} $ is 
\begin{align*} 
\mathcal{C}:= & \lbrace (\lambda, \Gamma, \Theta): \varepsilon\leq \lambda\leq M , \rho I\preceq \Theta -\Gamma\preceq I, 0\preceq\Gamma\preceq \alpha I, \\
& \Gamma\in\mathbf{D}_n,\Theta \in \mathbf{M}_n , \,\hat{\Sigma}^{-1}+\lambda^{-1}(\Theta-\Gamma)\succeq \beta I \rbrace 
\end{align*}
for a suitable $\beta >0$.

 Summing up we have the following 
\begin{teor}
Problem (\ref{Dual}) is equivalent to
\begin{equation}
\label{dual_compact}
\begin{aligned}
\min_{(\lambda,\Gamma, \Theta) \in \mathcal{C}} \tilde{J}(\lambda,\Gamma, \Theta).
\end{aligned}
\end{equation}
Both these problems admit solution.
\end{teor}

Before discussing uniqueness, it is convenient to further simplify the dual optimization problem: consider the function
\[F(\lambda,X):=-\lambda [\log(|\hat{\Sigma}^{-1}+\lambda^{-1}X|)
+\log|\hat{\Sigma}|-\delta]
\]
where $\lambda>0 $ and $X\in {\mathbf Q}_n$. Note that 
$$
F(\lambda,\Theta-\Gamma)=\tilde J (\lambda,\Gamma,\Theta).
$$
Moreover, $\Theta$ and $\Gamma$ are orthogonal over $\cal C$
 so that minimizing $\tilde J$ over ${\cal C}_0$ is equivalent to  minimize $F$
over the corresponding set
\begin{align*} 
\mathcal{C}_F:=  \lbrace (\lambda, X): & \  \lambda>0, \ 
X\in \mathbf{Q}_n, X\preceq I,  \\
& \ \chi(X)-X\succeq 0,  \ \hat{\Sigma}^{-1}+\lambda^{-1}X \succ 0 \rbrace.
\end{align*}
Therefore, from now on we can consider the following problem
\begin{equation}
\label{dual_simplif}
\begin{aligned}
\min_{(\lambda,X) \in \mathcal{C}_F} F(\lambda,X).
\end{aligned}
\end{equation}
Once obtained the optimal solution $(\lambda^*,X^*)$ we can recover the optimal values of the original multipliers simply by setting
$\Theta^*=\chi(X^*)$ and $\Gamma^*=\chi(X^*)-X^*$.

\subsection{Uniqueness of the solution of the dual problem}
The aim of this Section is to show that Problem (\ref{dual_simplif}) (and, hence
Problem (\ref{Dual})) admits a unique solution.
Since $\tilde J$ is the opposite of the dual objective function, $\tilde J$ is convex over $\cal C$.
It is then easy to check that $F$ is also a convex function over the convex set
$\mathcal{C}_F$. 
However, as we will see, $F$ is not strictly convex. Accordingly, establishing uniqueness of the minimum is not a trivial task.
 To this aim we need to compute the second variation $\delta^{2}F (\lambda,X;\delta \lambda, \delta X)$ of $F$ along all possible  directions
$(\delta \lambda, \delta X)$.
This second variation is a bilinear form that  can be represented in terms of an Hessian matrix of dimension $1+n^2$. It is possible to show that this matrix is positive definite and singular: its rank is $n^2$. This implies that $F$ is convex but there is a direction along which $F$ is not strictly convex.
It is then possible to prove that any optimal solution must be in the boundary of $\mathcal{C}_F$.
Moreover, if we consider the direction along which $F$ in non-strictly convex, the derivative of $F$ along this direction at any optimal point is non-zero.
In other words, there is a direction along which $F$ is affine but $F$ at any optimal point is not constant along this direction. Therefore we can prove the following result.

\begin{teor}
The dual problem admits a unique solution.
\end{teor}

\section{Recovering the solution of the primal problem} \label{sec_recovery}
Since the dual problem admits solution, we
 know that the duality gap between the primal and the dual problem is zero. The aim of this Section is to exploit this fact to recover the solution of the primal problem.\\
First, it is immediate to see that substituting the optimal solution of the dual problem $(\lambda^*,\Theta^*,\Gamma^*)$ into \eqref{der_var_wrt_sigma_opt} we obtain the optimal solution for $\Sigma$. The less trivial part is recovering the solution for $R$.\\
Zero duality gap allows us to apply the KKT theory and derive the following conditions:
\begin{equation}
\label{KKT_1}
\tr(\Lambda R)=0
\end{equation}
\begin{equation}
\label{KKT_2}
\tr(\Gamma(\Sigma-R))=0
\end{equation}
\begin{equation}
\label{KKT_3}
\tr(\Theta(\Sigma-R))=0.
\end{equation}

We start by considering \eqref{KKT_1}. It follows from \eqref{der_var_wrt_r} that
$\Lambda=I+\Gamma-\Theta $
where we notice that $\Lambda$ has no full rank. Therefore, we introduce its reduced singular value decomposition given by the following factorization
\begin{equation}
\label{reduced_svd}
\Lambda=U S U\tp 
\end{equation}
with $S\in \mathbb{R}^{n-r \times n-r}$, where $n-r$ is the rank of $\Lambda$. It follows that $U\tp U =I_{n-r}$. We plug \eqref{reduced_svd} in \eqref{KKT_1} and get
\begin{equation}
\begin{aligned}
0=\text{tr}[\Lambda R]= \text{tr}[USU\tp  R]
\Rightarrow U\tp RU=0.
\end{aligned}
\end{equation} 
 Therefore, by selecting a matrix $\tilde{U}$ whose columns form an orthonormal basis of $[{\rm im}(U)]^\bot$, we can
express  $R$ as:
\begin{equation}
\label{reduced_sdv_R}
R=\tilde{U}Q\tilde{U}\tp 
\end{equation}
with $Q=Q\tp \in \mathbb{R}^{r\times r}$. Notice that, the relationship $U\tp \tilde{U}=0$ holds, since the columns of $\tilde{U}$ form the orthogonal complement of the image of $U$.\\

By (\ref{KKT_3}), we know that $\Sigma-R$ is diagonal. Thus, we plug \eqref{reduced_sdv_R} into (\ref{KKT_3}) and obtain a linear system of equations: $\chi(\Sigma-\tilde{U}Q\tilde{U}\tp )=0$, or equivalently, 
\begin{equation}
\begin{aligned}
\chi(\tilde{U}Q\tilde{U}\tp )=\chi(\Sigma).
\end{aligned}
\end{equation}
In similar way, using (\ref{KKT_2}) we obtain an additional system of linear equations.
It is worth noting that the resulting system of equations always admits solution in $Q$ because the dual as well as the primal problem  admit solution. We conclude that, the linear system of equations admits a unique solution if and only if the solution of the primal which is unique.

\section{Conclusion}\label{sec_concl}
In this paper, factor analysis problem has been introduced for the realistic case in which the
covariance matrix of the data is estimated with an error which is not negligible. 
The dual analysis have been carried over to obtain a tractable mathematical problem that appears very promising for real applications.


\begin{thebibliography}{10}

\bibitem{agarwal2012noisy}
A.~Agarwal, S.~Negahban, and M.~J. Wainwright.
\newblock Noisy matrix decomposition via convex relaxation: Optimal rates in
  high dimensions.
\newblock {\em The Annals of Statistics}, pages 1171--1197, 2012.

\bibitem{ANDERSON_DEISTLER_1993}
B.~Anderson and M.~Deistler.
\newblock Identification of dynamic systems from noisy data: Single factor
  case.
\newblock {\em Mathematics of Control, Signals and Systems}, 6(1):10--29, 1993.

\bibitem{anderson1956statistical}
Theodore~W Anderson and Herman Rubin.
\newblock Statistical inference in factor analysis.
\newblock In {\em Proceedings of the third Berkeley symposium on mathematical
  statistics and probability}, volume~5, pages 111--150, 1956.

\bibitem{Bekker-deLeeuw_1987}
P.~A. Bekker and J.~de~Leeuw.
\newblock The rank of reduced dispersion matrices.
\newblock {\em Psychometrika}, 52(1):125?--135, 1987.

\bibitem{Bottegal-Picci}
G.~Bottegal and G.~Picci.
\newblock Modeling complex systems by generalized factor analysis.
\newblock {\em IEEE Transactions on Automatic Control}, 60(3):759--774, March
  2015.

\bibitem{BURT_1909}
C.~Burt.
\newblock Experimental tests of general intelligence.
\newblock {\em British Journal of Psychology, 1904-1920}, 3(1/2):94--177, 1909.

\bibitem{Chandrasekaran-Sanghavi-Parrilo-Willsky}
V.~Chandrasekaran, S.~Sanghavi, P.~A. Parrilo, and A.~S. Willsky.
\newblock Rank-sparsity incoherence for matrix decomposition.
\newblock {\em SIAM Journal on Optimization}, 21(2):572--596, 2011.

\bibitem{DEISTLER_2015}
M.~Deistler, W.~Scherer, and B.~Anderson.
\newblock The structure of generalized linear dynamic factor models.
\newblock In {\em Empirical Economic and Financial Research}, pages 379--400.
  Springer, 2015.

\bibitem{deistler2007}
M.~Deistler and C.~Zinner.
\newblock Modelling high-dimensional time series by generalized linear dynamic
  factor models: An introductory survey.
\newblock {\em Communications in Information \& Systems}, 7(2):153--166, 2007.

\bibitem{MIN_RANK_SHAPIRO_1982}
G.~Della~Riccia and A.~Shapiro.
\newblock Minimum rank and minimum trace of covariance matrices.
\newblock {\em Psychometrika}, 47:443--448, 1982.

\bibitem{ENGLE_ONE_FACTOR_1981}
R.~Engle and M.~Watson.
\newblock A one-factor multivariate time series model of metropolitan wage
  rates.
\newblock {\em Journal of the American Statistical Association},
  76(376):774--781, 1981.

\bibitem{FAZEL_MIN_RANK_APPLICATIONS_2002}
M.~Fazel.
\newblock Matrix rank minimization with applications.
\newblock {\em Elec. Eng. Dept. Stanford University}, 54:1--130, 2002.

\bibitem{Ferrante-Pav-Zorzi-12}
A.~Ferrante, M.~Pavon, and M.~Zorzi.
\newblock A maximum entropy enhancement for a family of high-resolution
  spectral estimators.
\newblock {\em IEEE Trans. Autom. Control}, 57(2):318--329, Feb. 2012.

\bibitem{forni_lippi_2001}
M.~Forni and M.~Lippi.
\newblock The generalized dynamic factor model: representation theory.
\newblock {\em Econometric theory}, 17(06):1113--1141, 2001.

\bibitem{GEWEKE_DYNAMIC}
J.~Geweke.
\newblock The dynamic factor analysis of economic time series models.
\newblock In {\em Latent Variables in Socio-Economic Models}, SSRI workshop
  series, pages 365--383. North-Holland, 1977.

\bibitem{ROBUSTNESS_HANSENSARGENT_2008}
L.~Hansen and T.~Sargent.
\newblock {\em Robustness}.
\newblock Princeton University Press, Princeton, NJ, 2008.

\bibitem{DEISTLER_1997}
C.~Heij, W.~Scherrer, and M.~Deistler.
\newblock System identification by dynamic factor models.
\newblock {\em SIAM Journal on Control and Optimization}, 35(6):1924--1951,
  1997.

\bibitem{ON_THE_PENA_BOX_MODEL_2001}
Y.~Hu and R.~Chou.
\newblock On the pena-box model.
\newblock {\em Journal of Time Series Analysis}, 25(6):811--830, 2004.

\bibitem{KALMAN_SELECTION_ECONOMETRICS_1983}
R.~Kalman.
\newblock {\em Identifiability and problems of model selection in
  econometrics}.
\newblock Cambridge University Press, 1983.

\bibitem{ROBUST_FILTERING}
B.~C. Levy and R.~Nikoukhah.
\newblock Robust state space filtering under incremental model perturbations
  subject to a relative entropy tolerance.
\newblock {\em IEEE Transactions on Automatic Control}, 58(3):682--695, 2013.

\bibitem{Pena_BOX__1987}
D.~Pena and G.~Box.
\newblock Identifying a simplifying in time series.
\newblock {\em Journal of the American Statistical Association},
  82(399):836--843, Sept. 1987.

\bibitem{PICCI_1987}
G.~Picci.
\newblock Parametrization of factor analysis models.
\newblock {\em Journal of Econometrics}, 41(1):17--38, 1989.

\bibitem{PICCI_PINZONI_1986}
G.~Picci and S.~Pinzoni.
\newblock Dynamic factor-analysis models for stationary processes.
\newblock {\em IMA Journal of Mathematical Control and Information},
  3(2-3):185--210, 1986.

\bibitem{SARGENT_SIM_1977}
T.~Sargent and C.~Sims.
\newblock {Business cycle modeling without pretending to have too much a priori
  economic theory}.
\newblock Technical Report~55, Federal Reserve Bank of Minneapolis, 1977.

\bibitem{scherrer1998structure}
Wolfgang Scherrer and Manfred Deistler.
\newblock A structure theory for linear dynamic errors-in-variables models.
\newblock {\em SIAM Journal on Control and Optimization}, 36(6):2148--2175,
  1998.

\bibitem{SCHUPPEN_1986}
J.~Schuppen.
\newblock Stochastic realization problems motivated by econometric modeling.
\newblock In C.~Byrnes and A.~Lindquist, editors, {\em Modeling Identification
  and Robust Control}, pages 259--275. North-Holland, 1986.

\bibitem{spearman_1904}
C.~Spearman.
\newblock {"General Intelligence," Objectively Determined and Measured}.
\newblock {\em American Journal of Psychology}, 15:201--293, 1904.

\bibitem{THURSTONE_1934}
L.~Thurstone.
\newblock The vectors of the mind.
\newblock {\em Psychological Review}, 41:1--12, 1934.

\bibitem{WATSON_ALGORITHMS_1983}
M.~Watson and R.~Engle.
\newblock {Alternative algorithms for the estimation of dynamic factor, mimic
  and varying coefficient regression models}.
\newblock {\em Journal of Econometrics}, 23(3):385--400, Dec. 1983.

\bibitem{ROBUST_TAU}
M.~Zorzi.
\newblock On the robustness of the {B}ayes and {W}iener estimators under model
  uncertainty.
\newblock {\em Automatica}, 83:133 -- 140, 2017.

\bibitem{ROBUST_KALMAN_2017}
M.~Zorzi.
\newblock Robust {K}alman filtering under model perturbations.
\newblock {\em IEEE Transactions on Automatic Control}, 62(6):2902--2907, 2017.

\bibitem{convtau}
M.~Zorzi.
\newblock Convergence analysis of a family of robust {K}alman filters based on
  the contraction principle.
\newblock {\em SIAM Journal on Control and Optimization}, 2017 (accepted).

\bibitem{BSL}
M.~Zorzi and A.~Chiuso.
\newblock Sparse plus low rank network identification: A nonparametric
  approach.
\newblock {\em Automatica}, 76:355 -- 366, 2017.

\bibitem{Zorzi-Ferrante-Automatica-12}
M.~Zorzi and A.~Ferrante.
\newblock On the estimation of structured covariance matrices.
\newblock {\em Automatica}, 48(9):2145--2151, Sep. 2012.

\bibitem{ROBUST_CONV_2015}
M.~Zorzi and B.~C. Levy.
\newblock On the convergence of a risk sensitive like filter.
\newblock In {\em 54th IEEE Conference on Decision and Control (CDC)}, pages
  4990--4995, 2015.

\bibitem{MFA}
M.~Zorzi and R.~Sepulchre.
\newblock Factor analysis of moving average processes.
\newblock In {\em European Control Conference (ECC)}, pages 3579--3584, 2015.

\bibitem{LATENTG}
M.~Zorzi and R.~Sepulchre.
\newblock {A}{R} identification of latent-variable graphical models.
\newblock {\em IEEE Trans. on Automatic Control}, 61(9):2327--2340, 2016.

\end{thebibliography}

\end{document}